\documentclass[3p, times, 11pt, authoryear, english, review]{elsarticle}

\usepackage{graphicx}
\usepackage{amssymb}
\usepackage{amsmath}
\usepackage{tikz}
\usepackage{floatrow}
\usepackage{blindtext}
\usepackage{hyperref}
\usepackage{subcaption}
\usetikzlibrary{bayesnet}
\newcommand{\bx}{\textbf{x}}

\newcommand{\bu}{\textbf{u}}
\newcommand{\bh}{\textbf{h}}
\newcommand{\bz}{\textbf{z}}
\newcommand{\bmu}{\boldsymbol{\mu}}
\newcommand{\bsigma}{\boldsymbol{\sigma}}
\newcommand{\bSigma}{\boldsymbol{\Sigma}}

\usepackage{algorithm}  
\usepackage{algorithmicx}  
\usepackage{algpseudocode}
\usepackage{url}
\usepackage{multirow}
\usepackage{array}

\journal{Transportation Research Part C: Emerging Technologies} 

\begin{document}
\begin{frontmatter}

\title{Predictive and Prescriptive Performance of Bike-Sharing Demand Forecasts for Inventory Management} 
\author[dtu]{Daniele Gammelli\corref{corresp}}
\ead{daga@dtu.dk} 

\address[dtu]{Department of Technology, Management and Economics, Technical University of Denmark, Kgs. Lyngby, Denmark, 2800}

\author[tum]{Yihua Wang}
\ead{yihua.wang@tum.de} 
\address[tum]{Logistics and Supply Chain Management, School of Management, Technical University of Munich, 80333 Munich, Germany}

\author[gron,twente]{Dennis Prak}
\ead{d.r.j.prak@utwente.nl} 
\address[gron]{Department of Operations, University of Groningen, PO Box 800, 9700 AV Groningen, The Netherlands}
\address[twente]{Department Industrial Engineering and Business Information Systems, University of Twente, PO Box 217, 7500 AE Enschede, The Netherlands}

\author[dtu]{Filipe Rodrigues}
\ead{rodr@dtu.dk}

\author[tum,mdsi]{Stefan Minner}
\ead{stefan.minner@tum.de} 
\address[mdsi]{Munich Data Science Institute (MDSI), Technical University of Munich, 85748 Garching, Germany}

\author[dtu]{Francisco Camara Pereira}
\ead{camara@dtu.dk}

\cortext[corresp]{Corresponding author.}
\date{\today}

\begin{abstract}
Bike-sharing systems are a rapidly developing mode of transportation and provide an efficient alternative to passive, motorized personal mobility.
The asymmetric nature of bike demand causes the need for rebalancing bike stations, which is typically done during nighttime. 
To determine the optimal starting inventory level of a station for a given day, a User Dissatisfaction Function (UDF) models user pickups and returns as non-homogeneous Poisson processes with piece-wise linear rates.
In this paper, we devise a deep generative model directly applicable in the UDF by introducing a \emph{variational Poisson recurrent neural network} model (VP-RNN) to forecast future pickup and return rates.
We empirically evaluate our approach against both traditional and learning-based forecasting methods on real trip travel data from the city of New York, USA, and show how our model outperforms benchmarks in terms of system efficiency and demand satisfaction. 
By explicitly focusing on the combination of decision-making algorithms with learning-based forecasting methods, we highlight a number of shortcomings in literature.     
Crucially, we show how more accurate predictions do not necessarily translate into better inventory decisions.
By providing insights into the interplay between forecasts, model assumptions, and decisions, we point out that forecasts and decision models should be carefully evaluated and harmonized to optimally control shared mobility systems.
\end{abstract}


\end{frontmatter}

\section{Introduction} 

The value of bike-sharing programs as an urban mobility solution is increasingly recognized by several cities around the world.
They provide a flexible transport solution that easily connects to other modalities, and mitigates traffic congestion and air pollution.
Bike-sharing concepts provide a healthy, cost- and time-efficient alternative to passive, motorized transportation \citep{sohrabi2020real}.
Whereas the first pioneering experiments -- such as the White Bikes project in Amsterdam (1965) -- were completely unregulated, successful later implementations depended heavily on IT to prevent vandalism and theft \citep{demaio2009bike}.
More recently, this IT usage enabled the application of advanced operations research and data science methods to optimize strategic, tactical, and operational decisions.
There are currently over 2000 bike-sharing programs active world-wide, covering over 9 million bikes, a rapid growth compared with 2 million bikes in 2016 and 700,000 bikes in 2013 \citep{Richter2018, Map2021}.
The vast majority of these projects are station-based and consist of networks of fixed-location stations with physical bike slots \citep{sohrabi2020real,shaheen2010bikesharing}.
One of the major challenges of bike-sharing systems is the spatio-temporal nature of mobility demand, such that trip origins and destinations are asymmetrically distributed (e.g. reflecting commuting into a downtown in the morning and vice-versa in the evening), making the overall system imbalanced and sensitive to disturbances. 
To counteract this, the bikes in such networks are usually rebalanced during the night, when demand is low. This is called static rebalancing \citep{laporte2015shared, tian2020rebalancing}.

The main operational-level decision problems of bike-sharing systems are demand forecasting, inventory decision-making, and rebalancing.
These three problems are typically considered sequentially, with inventory targets being constraints for the rebalancing (routing) problem, and demand forecasts in turn serving as inputs to decide on these inventory targets.
The target inventory level at the beginning of a day results from minimizing the so-called ``user dissatisfaction'', a penalty cost arising from arriving customers that do not find an available bike, and returning customers that do not find an empty slot to return their bike.
Theoretical papers on bike-sharing rebalancing either assume given (ranges of) starting inventory levels or given functions for the cost resulting from these starting inventory levels.
Empirical applications typically use a User Dissatisfaction Function (UDF). The most widely used UDF was first proposed by \citep{Raviv&Kolka2013}, modelling pickups and returns as non-homogeneous Poisson processes. 
The rates of these processes are assumed to be piece-wise constant (e.g. per hour) and obtained by taking historical averages of the same day and hour \citep{OMahony2015,Schuijbroek2017,freund2019analytics}.

A vast literature stream on demand forecasting for bike-sharing systems has rapidly emerged during the last decade.
The topic attracts attention as it is an exemplary case where historical demand records together with data on explanatory variables are abundantly available.
Bike-sharing demand is known to be heavily dependent on temporal information (intra- and inter-day), but also on the weather \citep{eren2020review}, with several machine learning approaches being applied to model these relationships. 
However, the resulting forecasts are typically studied in isolation from inventory decisions, where authors judge the quality of their forecasts on standard accuracy metrics, such as MAE, (R)MSE, and $\mbox{R}^2$, but not on their eventual performance in the UDF.
In this work, we argue that predictive and prescriptive performance goals should be carefully aligned when designing new predictive models, so to understand the relations between different methods and avoid unconscious overfitting to practically irrelevant forecasting metrics. 

The contribution of this paper is threefold.
First, we propose a neural architecture capable of modelling the pickups and returns as Poisson processes, thus being directly applicable in the UDF.
Specifically, we propose a deep generative model whereby we represent the unknown Poisson rates as latent variables and where the time-dependent dynamics are captured by a Recurrent Neural Network (RNN). 
Second, we empirically evaluate our model against both traditional and learning-based approaches on real trip data from the 30 most active stations of New York Citi Bike, and show how our model outperforms benchmarks in terms of predictive and prescriptive performance.
Third, we study existing mismatches between forecasting accuracy and decision performance.
Specifically, the user dissatisfaction cost corresponding to a certain starting inventory level is a complex function of all pickup and return rates during the day and the hourly differences between them, thus creating a misalignment between the prediction and decision objectives.
We propose to measure the error in the daily cumulative difference between pickups and return rates, and find that this better predicts inventory performance than MAE, MSE, and $\mbox{R}^2$.

The remainder of this paper is structured as follows. 
We first summarize relevant research directions in Section \ref{sec:literature_review}.
We then introduce the main theoretical foundations and formally present the proposed approach in Section \ref{sec:methodology}.
Lastly, we discuss empirical results on real world trip travel data in Section \ref{sec:experiments}. Section \ref{sec:conclusion} concludes the paper. 

\section{Literature review} 
\label{sec:literature_review}
We review the literature that is relevant to our study in three main streams: static rebalancing, inventory modelling, and demand forecasting. We discuss these streams top-down, indicating how rebalancing problems depends on inventory models, that in turn depend on demand forecasts.
Finally, we touch upon the disconnect between demand forecasts and inventory models for bike-sharing, and its parallels in the wider demand forecasting and inventory control literature.
For a recent, general literature review on bike-sharing problems, we refer to \cite{Shui2020}. 

\subsection{Static rebalancing}
The vast majority of bike-sharing systems are statically rebalanced during nighttime, when demand is low and the impact of rebalancing is highest \citep{laporte2015shared}.
As a result, the majority of existing literature focuses on static rebalancing \citep{tian2020rebalancing}.
Examples of work on dynamic rebalancing during the day are \citet{angelopoulos2018incentivized}, \citet{caggiani2018modeling}, and \citet{warrington2019two}.
Most literature on static rebalancing assumes a given target inventory level for each station as input to the routing problem \citep[e.g.][]{chemla2013bike, dell2014bike, wang2021enhanced}.
\citet{erdougan2014static} define a range of allowable inventory levels for each station, leading to more cost-efficient routing.
This approach is also adopted by \citet{kadri2016branch} and \cite{Schuijbroek2017}. 

Instead of assuming a given set of allowable starting inventory levels in the rebalancing problem, a more integrated approach is to include a UDF. This function maps the starting inventory level to a penalty cost, taking into account the stochastic transactions that occur throughout the day.
This approach is used by e.g. \citet{raviv2013static}, \citet{szeto2016chemical}, and \citet{ho2017hybrid}, although these authors do not further describe the UDF. They solely assume its existence.
Contrarily, \citet{vogel2014hybrid}, \citet{frade2015bike}, and \citet{Datner2019} assume that demand for trips between pairs of stations is deterministic. In real life, demand is uncertain, which makes the previously mentioned UDF approaches more suitable inputs to the rebalancing problem. 

\subsection{Inventory modelling}
In rebalancing models, a number of studies approach the question of how to define either the range of allowable starting inventory levels, or the UDF. \citet{nair2011fleet} and \citet{nair2013large} introduce two-sided constraints for failed pickups at empty stations and failed returns at full stations.
However, they only consider net demand and not the evolution of pickups and returns throughout the day. Also \citet{Maggioni2019} define a penalty based on net demand, which is assumed to be uniform, exponential, Gaussian, or log-normal.

\citet{Raviv&Kolka2013} model pickups and returns as independent, non-homogeneous Poisson processes with piece-wise constant rates.
\cite{Schuijbroek2017} use the same modelling approach, but combine it with the logic of \citet{nair2011fleet} to find a range of inventory levels that satisfy constraints with respect to the number of failed pickups and returns. Variants of this approach are commonly used to model the cost incurred at stations for given starting inventory levels in various applications.
Using a UDF with independent, non-homogeneous Poisson pickup and return processes, \citet{OMahony2015} determine rebalancing decisions, \citet{ccelebi2018bicycle} find the best locations for bike-sharing stations, and \citet{freund2019analytics} determine slot allocations and devise incentives for crowdsourcing rebalancing. 
Also in large-scale applications, where inventory levels and routing decisions are taken for multiple stations simultaneously, this UDF approach is used to model pickups and returns \citep[e.g.][]{jian2016simulation, alvarez2016optimizing}.

Almost all authors that empirically apply the UDF, use some form of historical averaging to estimate the pickup and return rates. \citet{Raviv&Kolka2013} use time intervals of 1, 5 and 30 minutes.
\citet{alvarez2016optimizing} and \citet{ccelebi2018bicycle} use an hourly interval, whereas \citet{Schuijbroek2017} and \citet{jian2016simulation} define the time interval as 15 and 30 minutes, respectively. \citet{OMahony2015} and \citet{freund2019analytics} use 20-minute intervals. 

In line with the reviewed literature, we use the UDF approach of \citet{Raviv&Kolka2013} to set and test a target inventory level.
For demonstration purposes, we deem comparing single target inventory levels more illustrative than comparing ranges of inventory levels and their implied service levels. Nevertheless, we remark that the procedure can be applied analogously to the service level approach of \cite{Schuijbroek2017}. 

\subsection{Demand forecasting}
Whereas empirical contributions to bike-sharing rebalancing and inventory decision making mainly use simple historical averaging to estimate the pickup and return rates, there exists a vast literature on bike-sharing demand forecasting. Attempts are made to derive explanatory power from exogenous variables, such as weather and temporal information, using both classical and machine learning prediction techniques.

\cite{rixey2013station} uses multivariate regression with data gathered from multiple bike-sharing systems, identifying a number of variables that have statistically significant correlations with station-level demand. More recently, the focus has shifted to machine learning approaches. System-level demand is forecasted by \cite{Xu2018} using long short-term memory neural networks, and by \cite{Guo2019} using graph neural networks. 

More closely resembling the rebalancing decisions that are to be made, several authors have also applied machine learning to forecast station-level demand. \cite{wang2018short} and \cite{chen2020predicting} employ RNNs, \cite{Lin2018} propose graph neural networks, and \cite{sohrabi2020real} use a generalized extreme value model. \cite{fournier2017sinusoidal} use a sinusoidal model to deal with seasonalities. Random forests have been adopted by \cite{yang2016mobility}, \cite{Du2019}, and \cite{ve2020season}. \cite{Gammelli2020, Gammelli2020EtAl} employ probabilistic techniques to estimate true demand using Tobit regression combined with Gaussian processes to mitigate the bias caused by censored demand observations, in both single and multi-output settings. In the presence of demand censoring, \cite{Negahban2019} estimates real demand with a combination of simulation and bootstrapping, whereas \cite{albinskiminner} present a data-driven approach to estimate achieved service levels.  \cite{boufidis2020development} compare various machine learning models in predicting station-level hourly pickups and returns. \cite{zhang2021data} iteratively update demand forecasts using a neural network and optimize a static rebalancing problem, but also separate both tasks.

Whereas several above-mentioned authors do focus on predicting station-level pickups and returns, they solely judge their forecasts using traditional loss measures for separate pickup and return demands.  It is not explored how (pickup and return) demand forecasts perform if they are used to optimize stations' starting inventory levels. This existing disconnect between prediction and optimization has been addressed in general terms by \citet{elmachtoub2021smart}. Particularly, the lacking interface between demand forecasting and inventory control has been pointed out by \citet{tratar2010joint}, \citet{prak2017calculation}, and \citet{kourentzes2020optimising}.
\citet{syntetos2010forecasting} state that the orders of magnitude of forecasting accuracy and inventory performance may differ wildly. \citet{babai2014intermittent} find that positively biased forecasts can actually be beneficial if the demand distribution is misspecified. 

In the field of bike-sharing, where (forecasts of) pickups and returns together determine the inventory trajectory of a station, the interface between predictive and prescriptive performance remains unstudied, despite the abundance of forecasting methods applied. This paper sheds light on how predictions affect decisions and, ultimately, system performance.

\section{Methodology}
\label{sec:methodology}
\begin{figure}[t]
\centering
\includegraphics[width=0.88\textwidth]{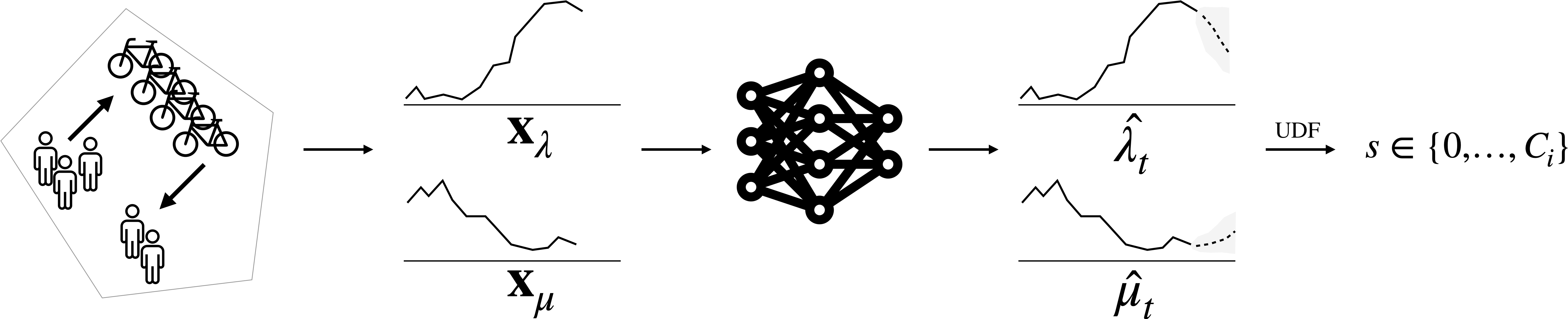}
\caption{An illustration of the framework determining the proposed inventory decision strategy.
Given historical observations of pickup ($\bx_{\mu}$) and return ($\bx_{\lambda}$) processes, the VP-RNN computes predictions ($\hat{\mu}_t$ and $\hat{\lambda}_t$). With these, we estimate and minimize the UDF to find the optimal daily starting inventory level ($s$).}
\label{fig:framework}
\end{figure}

In this section, we introduce a framework\footnote{Code available at: \url{https://github.com/DanieleGammelli/variational-poisson-rnn}} for inventory decision-making in bike-sharing systems.
As illustrated in Figure \ref{fig:framework}, this framework consists of a novel probabilistic neural architecture to estimate future pickup and return rates.
These estimates are then used in an inventory optimization model which defines pickups and returns as independent, non-homogeneous Poisson processes, and calculates the expected penalty due to failed pickups and returns as a function (the UDF) of the starting inventory level. 
Lastly, we decide on the optimal inventory level by minimizing the UDF.

In this section, we introduce the theoretical fundamentals of our proposed approach.
Specifically, we will first review an inventory model that utilizes the predicted pickup and return rates to determine the target inventory level at the beginning of a day (Section \ref{subsec:inventory_model}).  
We then outline the proposed VP-RNN (Section \ref{subsec:vprnn}) by first reviewing the theory and notation describing latent variable models, recurrent neural networks and approximate inference (Section \ref{subsubsec:lvms}, \ref{subsubsec:approximate_inference}, \ref{subsubsec:rnns}), on which we will build to introduce the proposed generative model (Section \ref{subsubsec:architecture}). 

\subsection{Inventory decision model}
\label{subsec:inventory_model}
In this section, we present an inventory model that uses the predicted pickup and return rates to determine the target inventory level at the beginning of a day. 
The inventory model uses the User Dissatisfaction Function (UDF) proposed in \citet{Raviv&Kolka2013}. 
A bike station is modeled as a double-ended $M_t/ M_t /1/C$ queuing system, with the number of customers in the queue representing the number of bikes in the station.
The customer inter-arrival times (for bike returns) and service times (i.e., inter-arrival times for bike pickups) are assumed to be exponentially distributed with time-dependent rates $\lambda_t$ and $\mu_t$, respectively. Similar to \cite{Raviv&Kolka2013}, we assume that these rates are piece-wise constant.
The capacity of the system $C$ represents the capacity of the bike station, i.e., the total number of docks in the station.

We consider a static rebalancing problem where bikes are rebalanced overnight.
Our goal is to determine the target inventory level for rebalancing operations in order to minimize the number of unsatisfied users for both bike pickup and return during the day.
According to \cite{Raviv&Kolka2013}, the procedure of determining the target inventory level with predicted pickup rates $\hat{\mu}_t$ and return rates $\hat{\lambda}_t$ can be divided into two steps as follows.

\begin{enumerate}
    \item Calculate the transient probability $p(s,\delta,t) \equiv {\rm Pr}\, (S(t) = \delta \,| \, S(0) = s)$, which is the probability of the station being at inventory level $\delta \in \{0,...,C\} $ at time $t \geq 0$ given that its starting inventory at time $0$ was $s$.
    In a non-stationary queue, the transient probabilities are solutions to the Kolmogorov forward Equation \eqref{eq:trans_prob}, which can be solved efficiently using the fourth-order Runge-Kutta method \citep{ross2014introduction}.
    \begin{align}
    \label{eq:trans_prob}
        \dot p(s,0,t) &= \hat{\mu}_t \cdot p(s,1,t) - \hat{\lambda}_t \cdot p(s,0,t) \nonumber\\
        \dot p(s,\sigma,t) &= \hat{\mu}_t \cdot p(s,\sigma+1,t) + \hat{\lambda}_t \cdot p(s,\sigma-1,t) - (\hat{\mu}_t+\hat{\lambda}_t) \cdot p(s,\sigma,t) && \sigma= 1, ..., C-1\\
        \dot p(s,C,t) &= \hat{\lambda}_t \cdot p(s,C-1,t) - \hat{\mu}_t \cdot p(s,C,t) \nonumber
    \end{align}
    \item Calculate the expected penalty (UDF) due to failed pickups and returns over the observation period $[0,T]$ for all possible starting inventories $s \in \{0,...,C\}$. The optimal starting inventory $s^*$ minimizes the UDF as shown in  \eqref{eq:inventory_decision}.
    \begin{align}  
    \label{eq:inventory_decision}
        UDF(s) & =\int_{0}^{T}l_p \cdot \hat{\mu}_t p(s,0,t) + l_r \cdot \hat{\lambda}_t p(s,C,t)\, {\rm d}t \nonumber\\
        s^* & = \mathop{\arg\min}\limits_{s} \, UDF(s)
    \end{align}
    Here, $l_p$ and $l_r$ denote the unit penalty for each lost pickup and lost return, respectively. The first term in the integral represents the expected user dissatisfaction accumulated when the station is empty, and the second term represents the expected user dissatisfaction accumulated when the station is full. 
\end{enumerate}

\subsection{Variational Poisson RNN}
\label{subsec:vprnn}
In this section, we first review and summarize key concepts on latent variable models (Section \ref{subsubsec:lvms}), approximate inference (Section \ref{subsubsec:approximate_inference}), and recurrent neural networks (Section \ref{subsubsec:rnns}).
We then build on these concepts to introduce the neural architecture for the proposed VP-RNN (Section \ref{subsubsec:architecture}).

\subsubsection{Latent Variable Models (LVMs)}
\label{subsubsec:lvms}
One of the central problems in the statistical sciences and machine learning is that of density estimation, i.e., the construction of a model of a probability distribution $p(\bx)$ given a finite sample of $N$ data points $\mathcal{D}: \{\bx_1, \ldots, \bx_N\}$ drawn from that distribution.
A traditional approach to the problem of density estimation involves a parametric model $p_{\theta}(\bx)$, in which a specific form for the density is proposed which contains a set of learnable parameters $\theta$. 
The parametric model of interest will be a Poisson distribution given by 
\begin{equation}
p_{\theta}(\bx) = \text{Pois}(\bx \, | \, \mathbf{\lambda}),
\end{equation}
where $\theta : \{\mathbf{\lambda} \}$ is the set of learnable parameters containing the rate of the Poisson distribution.
Learning, or parameter estimation, is then achieved by maximizing the (log) likelihood of the observed dataset as a function of the parameters, where it is assumed that the data points $\bx_i$ are drawn independently from $p(\bx)$.

Of particular interest for this paper is the concept of \emph{latent variables}. 
Specifically, rather than modelling $p(\bx)$ directly, we introduce a set of unobserved latent variables $\bz$ by expressing a model for the joint probability distribution $p(\bx, \bz)$. 
In practice, this is done by defining the joint probability as a product of two densities: the \emph{prior distribution} $p(\bz)$ and the \emph{likelihood} $p(\bx \, | \, \bz)$ (sometimes referred to as the \emph{sampling} or \emph{data distribution}), $p(\bx, \bz) = p(\bx \, | \, \bz)p(\bz)$.
In this context, parameter estimation, or \emph{inference}, is achieved by using Bayes' rule, yielding the following \emph{posterior} density:
\begin{equation}
\label{eq:posterior}
p(\bz \, | \, \bx) = \frac{p(\bx, \bz)}{p(\bx)} = \frac{p(\bx \, | \, \bz)p(\bz)}{p(\bx)},
\end{equation}
where $p(\bx) = \int p(\bx \, | \, \bz)p(\bz) d\bz$, and the integral is over all possible values of $\bz$ (or $p(\bx) = \sum_{\bz} p(\bx \, | \, \bz)p(\bz)$ in case of discrete $\bz$).

\subsubsection{Approximate Inference in LVMs}
\label{subsubsec:approximate_inference}
The posterior distribution in \eqref{eq:posterior} compactly represents our beliefs about the latent variables after having observed the data $\mathcal{D}$, and is a key component for probabilistic reasoning in LVMs.
In many cases of practical interest however, the posterior is intractable.
Specifically, this intractability often derives from the lack of an analytical solution for the integral appearing in the denominator of \eqref{eq:posterior}.
To address this intractability, we focus on deterministic techniques such as \emph{variational inference} (VI) \citep{JordanEtAl1999, Blei_2017, zhang2018advances}. At a high-level, in VI we use ideas from the calculus of variations to find a parametric approximation $q(\bz)$ that minimizes a measure of dissimilarity between $q(\bz)$ and the true, intractable posterior $p(\bz \, | \, \bx)$. 
Out of the many different ways to measure dissimilarity between two distributions, variational inference uses the \emph{Kullback-Leibler (KL)} divergence. That is, we are interested in minimizing the following divergence between the variational (or approximate) distribution $q(\bz)$ and the posterior distribution $p(\bz \, | \, \bx)$, defined as:
\begin{equation}
    \mathbb{KL}\left[q(\bz) || p(\bz \, | \, \bx)\right] = -\mathbb{E}_{q(\bz)} \left[ \log \frac{p(\bz \, | \, \bx)}{q(\bz)}\right],
    \label{eq:kl}
\end{equation}
where $\mathbb{E}_{q(\bz)}$ denotes an expectation over $q(\bz)$.
In order to define a tractable objective for our inference problem (i.e., one where the intractable posterior $p(\bz \, | \, \bx)$ does not appear in the formulation), we can rewrite \eqref{eq:kl} using \eqref{eq:posterior} (as well as the properties of the logarithm) as
\begin{align}
    \mathbb{KL}\left[q(\bz) || p(\bz \, | \, \bx)\right] &= -\mathbb{E}_{q(\bz)} \left[ \log \frac{p(\bx, \bz)}{q(\bz)} - \log p(\bx) \right] \\
    & = -\underbrace{\mathbb{E}_{q(\bz)} \left[\log \frac{p(\bx, \bz)}{q(\bz)}\right]}_{\mathcal{L}(q)} + \log p(\bx), 
\end{align}
where the marginal log-likelihood $\log p(\bx)$ can be taken out of the expectation because of its independence from $\bz$.
The quantity $\mathcal{L}(q)$ is know as \emph{Evidence Lower Bound} (ELBO) and represents a lower bound on the marginal log-likelihood, or evidence, $\log p(\bx)$, i.e. $\log p(\bx) \geq \mathcal{L}(q)$ for all $q(\bz)$. 
Concretely, this reformulation gives us a way to minimize the $\mathbb{KL}\left[q(\bz) || p(\bz \, | \, \bx)\right]$ by maximizing the ELBO with respect to the distribution $q(\bz)$, and therefore find the variational distribution best approximating the unknown posterior.
In other words, the closer the ELBO is to the marginal log-likelihood, the closer (in KL sense) the variational approximation will be to the posterior distribution.
Thus, variational methods allow us to reduce an inference problem into an optimization problem.

In practice, the variational distribution $q(\bz)$ is often restricted to a known parametric family for which the ELBO is tractable or simple to approximate, such as a Gaussian distribution. 
Thus, the maximization of the ELBO refers to a maximization with respect to the parameters $\phi$ of the variational distribution $q_{\phi}$ (e.g. $q_{\phi}(\bz) = \mathcal{N}(\bz \,|\, \phi)$, where $\phi = \{\bmu, \bSigma\}$ in the case of a Gaussian approximation).

In traditional variational inference, we learn a distinct set of parameters $\phi_i$ for each data point $\{\bx_i\}_{i=1}^N$, which can be problematic when facing large, high-dimensional datasets.
To avoid the linear growth in parameters with the number of data points, \emph{amortized inference} offers a viable alternative.
Specifically, rather than defining a set of parameters $\phi_i$ for each data point, amortized inference shares a unique set of parameters $\phi$ across all data points - thus, \emph{amortizing} the cost of variational inference.
As in the case of Variational Autoencoders (VAE) \citep{KingmaEtAl2014, RezendeEtAl2014}, we define an inference network, also known as encoder, that allows us to compute the parameters of the posterior approximation for any given data point. 
Specifically, in the case of a (diagonal) Gaussian variational approximation, we define an inference network with output characterizing the mean and variance vectors as:
\begin{align}
    q_{\phi}(\bz_i \, | \, \bx) & = \mathcal{N}(\bz_i \, | \, \bmu_i, \bsigma_i^2 I), \\
    [\bmu_i, \bsigma_i^2] & = f_{\phi}(\bx_i), \nonumber
\end{align}
where $f_{\phi}$ can be any parametric function such as a deep neural network, and $I$ is the identity matrix.

\subsubsection{Recurrent Neural Networks}
\label{subsubsec:rnns}
We summarize the usage of recurrent neural networks (RNNs) for sequential data modelling. RNNs are widely used to model variable-length sequences $\mathbf{x} = (\mathbf{x}_1, \mathbf{x}_2, \ldots, \mathbf{x}_T)$, possibly influenced by external covariates $\mathbf{u} = (\mathbf{u}_1, \mathbf{u}_2, \ldots, \mathbf{u}_T)$. 
The core assumption underlying these models is that all observations $\mathbf{x}_{1:t}$ up to time $t$ can be summarized by a learned deterministic representation $\mathbf{h}_t$. 
At any timestep $t$, an RNN recursively updates its hidden state $\mathbf{h}_t \in \mathbb{R}^p$ by computing:
\begin{equation}
    \mathbf{h}_t = f_{\theta_{\bh}}(\mathbf{u}_t, \mathbf{h}_{t-1}),
\end{equation}
where $f$ is a deterministic non-linear transition function parametrized by $\theta_{\bh}$, such as an Long Short-Term Memory (LSTM) cell or a Gated Recurrent Unit (GRU).
The sequence is then modelled by defining a factorization of the joint probability distribution as the following product of conditional probabilities:
\begin{align}
    p(\mathbf{x}_1, \mathbf{x}_2, \ldots \mathbf{x}_T) &= \prod_{t=1}^{T}{p(\mathbf{x}_t|\mathbf{x}_{<t})} \nonumber\\
    p(\mathbf{x}_t|\mathbf{x}_{<t}) &= g_{\theta_{\bx}}(\mathbf{h}_t), \label{eq:rnn}
\end{align}
where $g$ is typically a non-linear function with parameters $\theta_{\bx}$.

\subsubsection{VP-RNN neural architecture}
\label{subsubsec:architecture}

\begin{figure*}[t]
\begin{subfigure}{0.45\textwidth}
\centering
\includegraphics[width=0.65\textwidth]{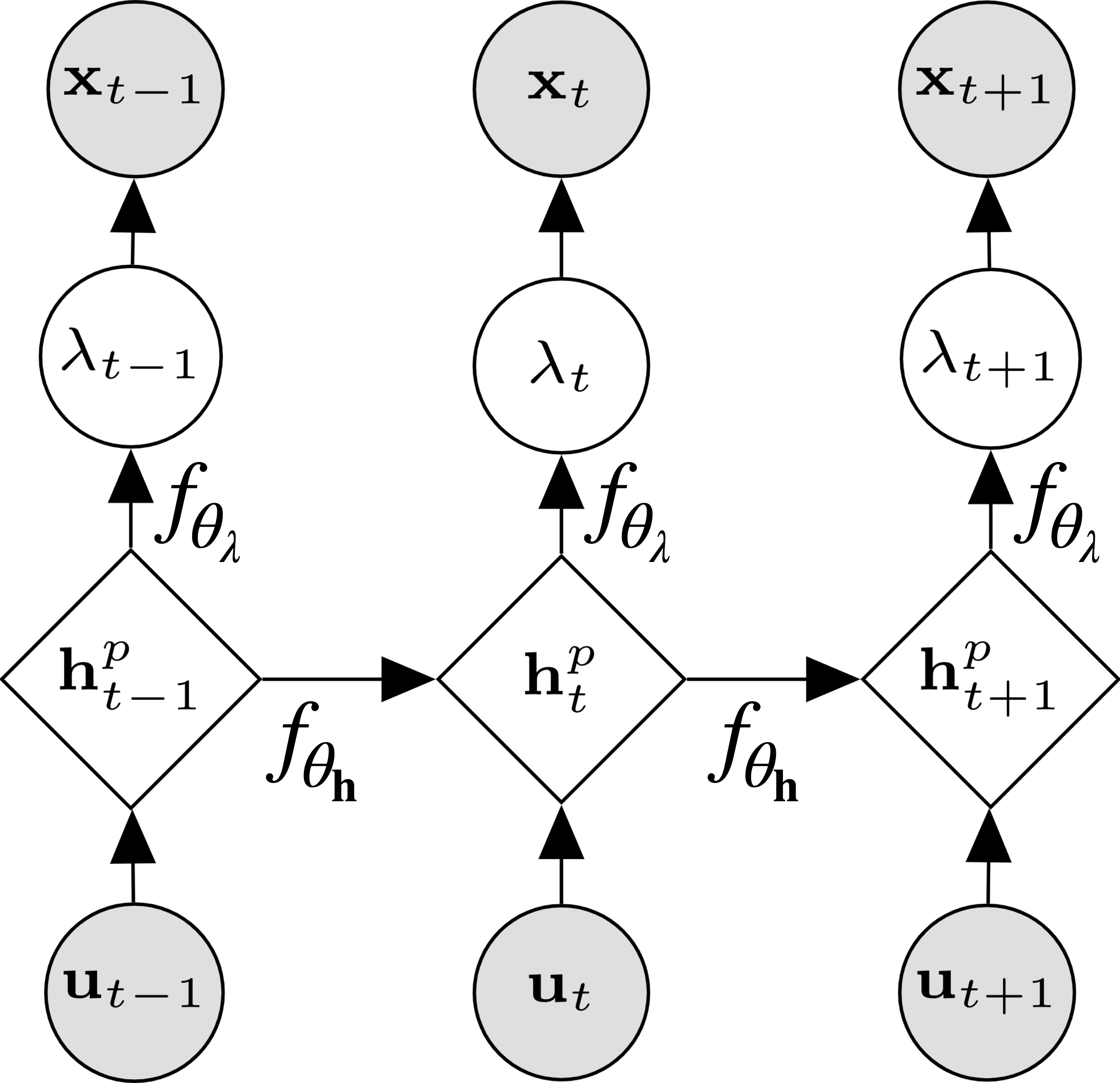}

\caption{Generative model}
\end{subfigure}\hspace{\fill}
~
\begin{subfigure}{0.45\textwidth}
\centering
\includegraphics[width=0.65\textwidth]{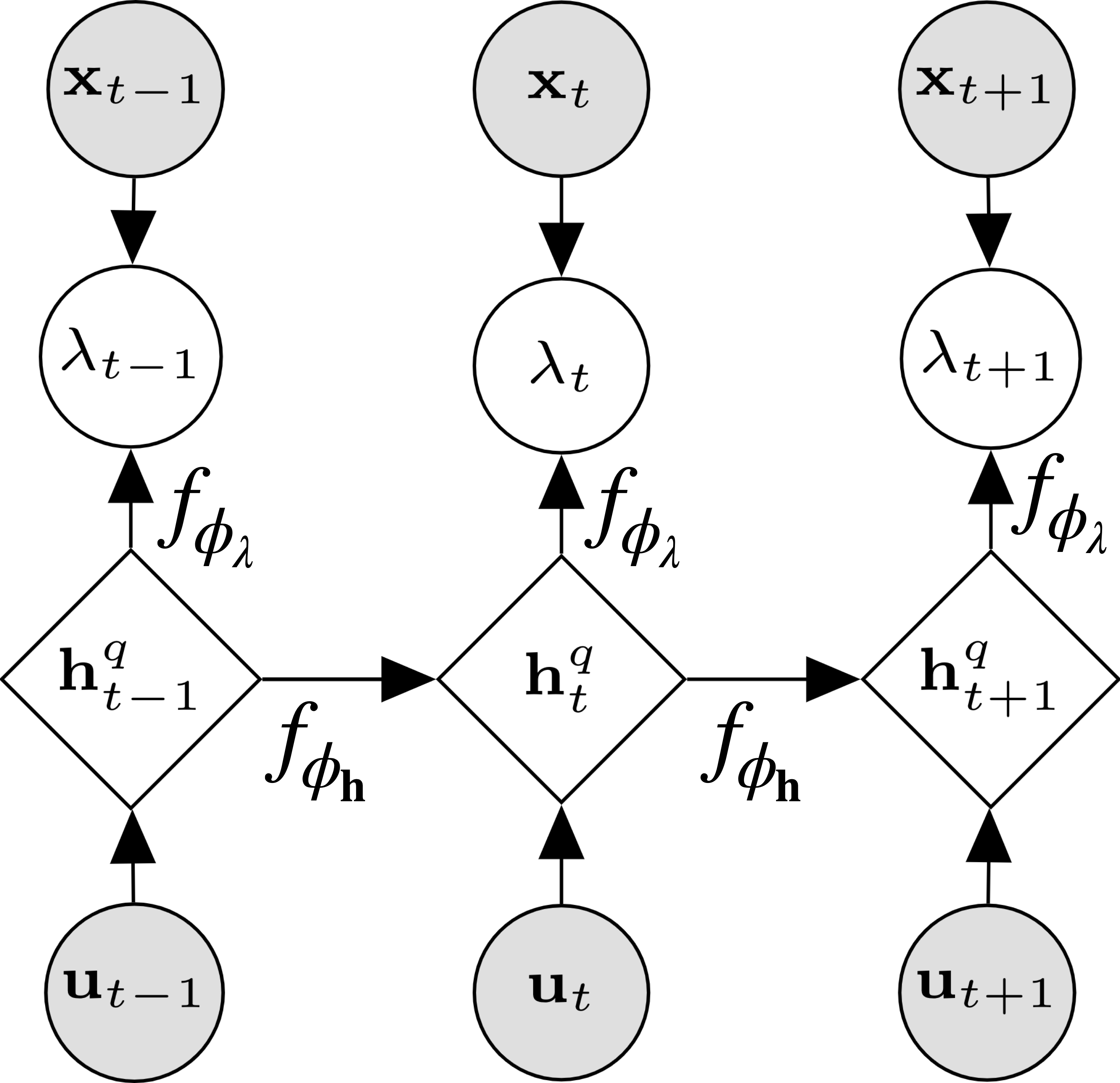}
\caption{Inference Network}
\end{subfigure}\hspace{\fill}
\caption{Graphical representation of the generative model $p_{\theta}$ (left) and inference network $q_{\phi}$ (right) characterizing the proposed VP-RNN. 
Shaded nodes represent observed variables, while non-shaded nodes represent either deterministic (diamond-shaped) or latent (circle-shaped) variables.}
\label{fig:pgm}
\end{figure*}

In this section, we define the generative model $p_{\theta}$ and inference network $q_{\phi}$ characterizing the proposed Variational Poisson RNN (VP-RNN) for the purpose of pickup and return demand modelling.
A schematic illustration of the VP-RNN is shown in Figure \ref{fig:pgm}.

\smallskip \noindent\textbf{Generative model}: We assume data $\mathbf{X} = (\mathbf{x}_1, \mathbf{x}_2, \ldots, \mathbf{x}_T)$ to represent a sequence of realizations of a Poisson process over the pick-up (or, the return) process. That is, we recognize that our data is represented by count variables taking values $\mathbf{x}_t \in \mathbb{N} \cup \{0\}$ that we wish to model using a Poisson process specified by an appropriate rate parameter $\lambda_t$.
For simplicity, we will always refer to a generic rate parameter $\boldsymbol{\lambda}_t$ when describing the proposed VP-RNN. 
However, it is important to underline how $\boldsymbol{\lambda}_t$ can represent any arbitrary Poisson rate in the context of bike-sharing demand prediction, such as independent pickup and return rates $\mu_t$, $\lambda_t$, or even a 2-dimensional rate $\boldsymbol{\lambda}_t = [\mu_t, \lambda_t]$, jointly modeling the pickup and return processes.
We represent the rate $\boldsymbol{\lambda}_t$ as a latent variable whose time-dependent dynamics are modelled through an RNN. 
Specifically, the VP-RNN defines the following factorization of the joint probability distribution:
\begin{align}
    p(\bx_{1:T}, \boldsymbol{\lambda}_{1:T}, \bh^p_{1:T} \, | \, \bu_{1:T}, \bh^p_0) & = \prod_{t=1}^T p(\bx_t \, | \, \boldsymbol{\lambda}_t) \, p_{\theta_{\lambda}}(\boldsymbol{\lambda}_t \, | \, \bh^p_{t}) \,p_{\theta_{\bh}}(\bh^p_t \, | \, \bh^p_{t-1}, \bu_t), \\
    p(\bx_t \, | \, \boldsymbol{\lambda}_t) & = \text{Pois}(\bx_t \, | \, \boldsymbol{\lambda}_t) \nonumber \\
    p_{\theta_{\lambda}}(\boldsymbol{\lambda}_t \, | \, \bh^p_{t}) & = \mathcal{N}\left(\boldsymbol{\lambda}_t \, | \, \bmu_{0,t} \text{diag}(\bsigma_{0,t}^2)\right), \text{with   } [\bmu_{0,t}, \bsigma_{0,t}] = f_{\theta_{\lambda}}(\bh^p_{t}), \nonumber
\end{align}
where $\bmu_{0,t}$ and $\bsigma_{0,t}$ represent the parameters of the conditional prior distribution over the latent variable $\boldsymbol{\lambda}_t$ and where we assume $p_{\theta_{\bh}}(\bh^p_t \,|\, \bh^p_{t-1}, \bu_t) = \delta(\bh^p_t - \tilde{\bh})$, i.e. $\bh^p_t$ follows a delta distribution centered in $\tilde{\bh}_t = f_{\theta_{\bh}}(\bh^p_{t-1}, \bu_t)$.
In our implementation, $f_{\theta_{\lambda}}$ and $f_{\theta_{\bh}}$ are respectively a feed-forward neural network and a GRU cell with parameters $\theta_{\lambda}$ and $\theta_{\bh}$.
At its core, the VP-RNN exploits the representational power of RNNs to capture potentially complex long-term dependencies in the temporal evolution of demand. It then leverages the learned representation $\bh^p_t$ as conditioning variable for the conditional prior distribution over the Poisson rate variable $\boldsymbol{\lambda}_t$. 

\smallskip \noindent\textbf{Inference}: 
The variational approximation defining the VP-RNN directly follows the generative model's factorization as follows:
\begin{align}
    q_{\phi}(\boldsymbol{\lambda}_{1:T} \, | \, \bx_{1:T}) & = \prod_{t=1}^T q_{\phi}(\boldsymbol{\lambda}_t \, | \, \bx_t, \bh_t^q, \bu_t), \label{eq:inference_net}\\
    q_{\phi}(\boldsymbol{\lambda}_t \, | \, \bx_t, \bh_t^q, \bu_t) & = q_{\phi_{\lambda}}(\boldsymbol{\lambda}_t \, | \,  \bh_t^q) \, q_{\phi_{\bh}}(\bh_t^q \, | \, \bh_{t-1}^q, \bu_t), \nonumber\\ 
    q_{\phi_{\lambda}}(\boldsymbol{\lambda}_t \, | \,  \bh_t^q) & = \mathcal{N}\left(\boldsymbol{\lambda}_t \, | \, \bmu_{\lambda,t}, \text{ diag}(\bsigma_{\lambda,t}^2)\right), \text{with   } [\bmu_{\lambda,t}, \bsigma_{\lambda,t}] = f_{\phi_{\lambda}}(\bh_{t}^q), \nonumber
\end{align}
where $q_{\phi_{\bh}}(\bh_t^q \,|\, \bh_{t-1}^q, \bu_t)$ follows a delta distribution centered in $\tilde{\bh}_t^q = f_{\phi_{\bh}}(\bh_{t-1}^q, \bu_t)$.
Concretely, $f_{\phi_{\lambda}}$ and $f_{\phi_{\bh}}$ together describe the encoder network defining the parameters $\bmu_{\lambda,t}$ and $\bsigma_{\lambda,t}$ of the approximate posterior distribution.
In our implementation, $f_{\phi_{\lambda}}$ and $f_{\phi_{\bh}}$ are respectively a feed-forward neural network and an LSTM cell with parameters $\phi_{\lambda}$ and $\phi_{\bh}$.
By explicitly resembling the model's factorization, the inference network defined in \eqref{eq:inference_net} also exhibits an implicit dependence on the entire history of $\bx_{1:t}$ and $\bu_{1:t}$ through $\bh_t^q$.
This implicit dependency on all information from the past can be considered as resembling a \emph{filtering} approach from the state-space model literature \citep{KoopmanEtAl2001}.
Denoting $\theta$ and $\phi$ as the set of model and variational parameters respectively, variational inference offers a scheme for jointly optimizing parameters $\theta$, $\phi$ and computing an approximation to the posterior distribution by maximizing the following step-wise evidence lower bound\footnote{A complete derivation is provided in the Appendix} (i.e. ELBO) through gradient ascent:
\begin{align}
    \mathcal{L}(\theta, \phi) & =  \mathbb{E}_{q_{\phi}(\boldsymbol{\lambda}_{1:T} \, | \, \bx_{1:T})} \left[\sum_{t=1}^{T} \log p_{\theta}(\bx_{t} \, | \, \boldsymbol{\lambda}_{t}) + \log p_{\theta}(\bh_{t}^p \, | \, \bh_{t-1}^p, \bu_{t}) \right] \label{eq:elbo}\\
    & - \sum_{t=1}^{T} \mathbb{KL} \left(q_{\phi}(\boldsymbol{\lambda}_{t} \, | \, \bh_{t}^q, \bx_{t}, \bu_t) \mid\mid p_{\theta}(\boldsymbol{\lambda}_{t} \, | \, \bh^p_{t}) \right). \nonumber
\end{align}

\section{Empirical Results}
\label{sec:experiments}
In this section, we demonstrate the performance of our proposed approach.
Specifically, the goal is to answer the following questions: (1) Can we learn to reliably predict future pickup and return rates? (2) Does predictive performance align with decision-making performance? (3) In case of a misalignment, what aspects should be taken into consideration when working on frameworks combining prediction and decision-making? To answer these questions, we first analyse the performance of the proposed VP-RNN in predicting pickup and return processes of bike-sharing demand against other learning-based approaches (Section \ref{subsec:bike_sharing_demand_prediction}). 
We then explicitly evaluate the predictions when used for inventory management tasks both quantitatively (Section \ref{subsec:inventory_management}) and qualitatively (Section \ref{subsec:qualitative_understanding}).

\smallskip We use a real-world dataset from New York Citi Bike \citep[][]{Citibike}. 
Citi Bike operates a station-based system, whereby the user of the service is not free to pick up or
drop off a bike in any location, but is restricted to a certain number of physical stations around New York. Our objective is to model the temporal evolution of station-level pickup and return demand in the bike-sharing system and use this understanding to decide on effective starting inventory levels. 

In all our experiments, we use data from the 30 most active stations in the Citi Bike's system from 1 January 2018 until 31 December 2018, covering approximately 25\% of all rides. 
As of December 2018, the system consisted of approximately 11,500 bikes, with 147,090 total annual memberships and an average demand in December 2018 of 41,172 rides per day.
The stations which we consider in this work are representative of a number of different demand patterns, such as morning pickup (return) peaks and evening return (pickup) peaks in e.g., residential (business) areas, as well as more balanced situations.

The data consists of individual records of users renting and returning bikes, which we aggregate to three distinct temporal aggregation levels: 15-, 30- and 60-minute intervals.
Once aggregated, the data at our disposal is characterized by the time series of station-level pickups and returns, which we aim to predict one-day ahead at the start of every new day, to reflect the decision that is to be made.
For all stations, we split the 12 months of data into train, validation and test sets using a ratio of 9/1/2 months, which we use respectively for training, model selection and early stopping, and the final evaluation of the implemented models.

For all models, we consider additional external explanatory variables to encode both meteorological and temporal information.
Specifically, our features are characterized by the following sources of information: (i) temperature [°C], (ii) probability of rain $\in [0,1]$, (iii) Day-of-Week (DoW), (iv) Time-of-Day (ToD), where we express both (iii) and (iv) as one-hot-encoded vectors.
We use hourly weather measurements as recorded by the National Climatic Data Center \citep{NCDC2018}.
In case of smaller temporal aggregation (15- and 30-minute intervals), we assume the weather measurements to remain constant throughout the hour.

\subsection{Predictive Results}
\label{subsec:bike_sharing_demand_prediction}
In this section, we analyze the performance of the proposed model on the task of pickup and return bike-sharing demand prediction.
We compare the performance of the proposed VP-RNN with other learning-based approaches that are often used in empirical bike-sharing literature. 
We further place this comparison in the context of an ablation study to better analyze the contribution of each individual component of the VP-RNN.
Concretely, we compare the performance of the following models:
\begin{enumerate}
    \item Historical Average (HA): given a temporal aggregation (e.g. 60-minute), the historical average for every combination of day-of-week and time-of-day (e.g. Monday 8 am) is calculated.
    
    \item Moving Average (MA): functionally equivalent to HA, with the only difference that the average is computed using only the last month of data in a rolling window.
    
    \item Linear Regression (LR): parametrizes the dependency of the number of pickups/returns $\bx_t$ on explanatory features $\bu_t$ through a linear relationship, estimated by ordinary least squares. 
    \item Poisson RNN (P-RNN): variation on the proposed VP-RNN not including an explicit latent variable over the Poisson rate $\boldsymbol{\lambda}_t$. In line with Section \ref{sec:methodology}, the P-RNN defines the following factorization of the joint distribution:
    \begin{align}
    p(\bx_{1:T}, \boldsymbol{\lambda}_{1:T}, \bh^p_{1:T} \, | \, \bu_{1:T}, \bh^p_0) & = \prod_{t=1}^T p(\bx_t \, | \, \boldsymbol{\lambda}_t) \, p_{\theta_{\lambda}}(\boldsymbol{\lambda}_t \, | \, \bh^p_{t}) \,p_{\theta_{\bh}}(\bh^p_t \, | \, \bh^p_{t-1}, \bu_t), \\
    p(\bx_t \, | \, \boldsymbol{\lambda}_t) & = \text{Pois}(\bx_t \, | \, \boldsymbol{\lambda}_t),
    \text{with   } \boldsymbol{\lambda}_t = f_{\theta_{\lambda}}(\bh^p_{t}). \nonumber
    \end{align}
    Given the absence of latent variables, the P-RNN allows for exact maximum likelihood estimation of the parameters.
    For all stations considered, in our implementation $p_{\theta_h}$ is a GRU with 128 hidden units, and $p_{\theta_{\lambda}}$ is a 2-layer MLP with 128 hidden units per hidden layer. 
    
    \item Variational Poisson RNN (VP-RNN): the model as described in Section \ref{subsubsec:architecture}.
    Similarly to P-RNN, in our implementation $p_{\theta_h}$ is a GRU with 128 hidden units, and $p_{\theta_{\lambda}}$ is a 2-layer MLP with 128 hidden units per hidden layer.
    However, the VP-RNN defines a Gaussian distribution over $\lambda$, opposed to a single point-estimate as in P-RNN. The inference network mirrors the implementation of the generative model where $q_{\theta_h}$ is a GRU with 128 hidden units, and $q_{\theta_{\lambda}}$ is a 2-layer MLP with 128 hidden units per hidden layer. 
    
    \item Multi-Output Variational Poisson RNN (MOVP-RNN): multi-output extension to the proposed VP-RNN. Specifically, this formulation allows to jointly model pickup and return processes by defining a multivariate regression variable $\bx_t = [x_{\mu, t}, x_{\lambda, t}]$, where $x_{\mu, t}$ and $x_{\lambda, t}$ represent pickup and return counts, respectively. By doing so, MOVP-RNN can potentially leverage correlations between the pickup and return temporal patterns. Our MOVP-RNN implementation uses the same number of model parameters as our VP-RNN.
\end{enumerate}

\begin{table}[t]
\small
\caption{Test prediction performance. We report average (std. dev.) performance over all stations considered.}
\tabcolsep=0.11cm
\begin{tabular}{c l c c c | c c c}
& & & Pickup & & & Return & \\
\hline
Aggregation & Models & RMSE & MAE & $\text{R}^2$ & RMSE & MAE & $\text{R}^2$\\
\hline
\multirow{6}{*}{60 min} & Historical Average      & 6.76 (2.40) & 4.19 (1.07) & 0.23 (0.33) & 6.65 (2.08) & 4.17 (0.99) & 0.29 (0.31) \\
& Moving Average          & 5.77 (2.00) & 3.45 (0.80) & 0.47 (0.07) & 5.80 (1.86) & 3.47 (0.75) & 0.50 (0.07) \\
& Linear Regression       & 6.67 (2.39) & 4.46 (1.35) & 0.25 (0.33) & 6.67 (2.20) & 4.50 (1.25) & 0.29 (0.31) \\
& Poisson-RNN             & 4.25 (1.15) & 2.65 (0.55) & 0.70 (0.05) & 4.28 (1.15) & 2.67 (0.51) & 0.71 (0.07) \\
& Variational Poisson-RNN & 3.91 (0.99) & 2.47 (0.46) & 0.74 (0.06) & 3.92 (0.93) & 2.49 (0.49) & 0.75 (0.07) \\
& Multi-Output VP-RNN     & \textbf{3.77 (0.96)} & \textbf{2.39 (0.46)} & \textbf{0.76 (0.06)} & \textbf{3.71 (0.82)} & \textbf{2.36 (0.41)} & \textbf{0.78 (0.06)} \\
\hline
\multirow{6}{*}{30 min} & Historical Average      & 3.61 (1.14) & 2.28 (0.49) & 0.25 (0.24) & 3.57 (1)    & 2.29 (0.46) & 0.3 (0.24)  \\
& Moving Average          & 3.28 (1.01) & 2.04 (0.42) & 0.4 (0.07)  & 3.3 (0.94)  & 2.05 (0.4)  & 0.43 (0.08) \\
& Linear Regression       & 3.7 (1.19)  & 2.5 (0.67)  & 0.21 (0.28) & 3.72 (1.13) & 2.52 (0.63) & 0.25 (0.26) \\
& Poisson-RNN             & 2.53 (0.53) & 1.61 (0.25) & 0.63 (0.08) & 2.53 (0.47) & 1.61 (0.24) & 0.64 (0.09) \\
& Variational Poisson-RNN & 2.39 (0.43) & 1.55 (0.22) & 0.66 (0.09) & 2.41 (0.42) & 1.56 (0.22) & 0.67 (0.09) \\
& Multi-Output VP-RNN     & \textbf{2.32 (0.42)} & \textbf{1.5 (0.21)} & \textbf{0.68 (0.09)} & \textbf{2.33 (0.4)} & \textbf{1.51 (0.21)} & \textbf{0.7 (0.08)} \\
\hline
\multirow{6}{*}{15 min} & Historical Average      & 2.07 (0.56) & 1.35 (0.25) & 0.21 (0.20) & 2.05 (0.50) & 1.35 (0.24) & 0.25 (0.20) \\
& Moving Average          & 1.94 (0.50) & 1.24 (0.23) & 0.32 (0.07) & 1.95 (0.47) & 1.25 (0.22) & 0.34 (0.09) \\
& Linear Regression       & 2.11 (0.59) & 1.44 (0.33) & 0.18 (0.23) & 2.12 (0.56) & 1.45 (0.32) & 0.21 (0.22) \\
& Poisson-RNN             & 1.59 (0.28) & 1.03 (0.15) & 0.52 (0.10) & 1.59 (0.26) & 1.03 (0.15) & 0.54 (0.11) \\
& Variational Poisson-RNN & 1.55 (0.25) & 1.02 (0.13) & 0.54 (0.11) & 1.56 (0.23) & 1.03 (0.13) & 0.56 (0.11) \\
& Multi-Output VP-RNN     & \textbf{1.53 (0.25)} & \textbf{1.00 (0.13)} & \textbf{0.55 (0.10)} & \textbf{1.53 (0.23)} & \textbf{1.01 (0.14)} & \textbf{0.57 (0.10)} \\
\hline
\end{tabular}
\label{tab:pred_results}
\end{table}

Table \ref{tab:pred_results} shows the predictive performance of the implemented models based on three commonly-used measures: Root Mean Squared Error (RMSE), Mean Absolute Error (MAE) and the coefficient of determination ($\text{R}^2$).
We now concentrate on the results for the 60-minute aggregation, as presented in Table \ref{tab:pred_results}, because they are representative also of the results for the other two temporal aggregations.
Unsurprisingly, results show how RNN-based approaches have a clear advantage when compared to the classical benchmarks that are typically used in empirical bike-sharing literature.
Table \ref{tab:pred_results} further highlights the contributions of each individual component of our proposed model.
First, results show how the MOVP-RNN is able to exploit its additional flexibility in modelling correlations between the pickup and return processes, obtaining better performance compared to its single-output variant (VP-RNN) across all metrics. 
Moreover, Table \ref{tab:pred_results} also highlights the gains in defining explicit latent variables over the rate parameter, thus allowing for a structured treatment of uncertainty and ultimately leading to more accurate predictions.

\begin{figure}[t]
\centering
\includegraphics[width=\textwidth]{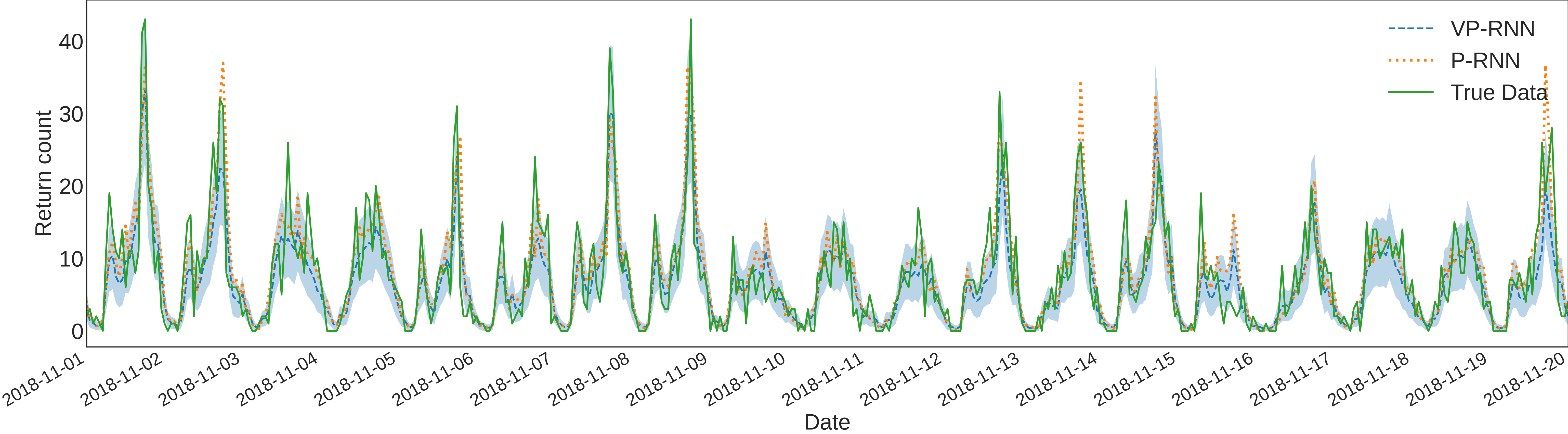}
\caption{A graphical representation of model prediction for 20 consecutive test days in station 504 of the Citi Bike system. The plot compares VP-RNN (blue, dashed curve), P-RNN (orange, dotted curve) against true realizations of bike return counts (green, continuous curve). The blue shaded area represents the 95\% interval under the posterior predictive distribution over the return Poisson rate parameter. The plot shows how VP-RNN is able to obtain higher log-likelihood values by averaging over multiple possible rates that could have generated the data.}
\label{fig:prnn_vs_vprnn}
\end{figure}

\begin{table}[t]
\small
\caption{Average test log-likelihood across stations. For the non-deterministic model (VP-RNN) the approximation of the marginal log-likelihood is indicated with the $\approx$ sign.}
\begin{tabular}{l c c }
\hline
Models & Pickup & Return\\
\hline
Poisson-RNN (P-RNN)              &  3666  & 4012 \\
Variational Poisson-RNN (VP-RNN) $\approx$ &  \textbf{3896}  & \textbf{4210}\\
\hline
\end{tabular}
\label{tab:pred_log_likelihood}
\end{table}

To further illustrate the potential advantages of explicitly modelling the rate as a latent variable, in Table \ref{tab:pred_log_likelihood} we compare the VP-RNN with its deterministic counterpart P-RNN.
Specifically, results show the test log-likelihood averaged over all stations of interest for the 60-minute aggregation.
We report exact log-likelihoods for P-RNN, while in the case of VP-RNN, we report the importance sampling approximation to the marginal log-likelihood using 30 samples, as in \cite{RezendeEtAl2014}.
For both pickup and return processes, we see how the combination of RNNs with latent variable models allows the VP-RNN to better estimate the demand process, obtaining higher log-likelihood values on held-out data.

Crucially, by explicitly allowing for the presence of latent variables, the VP-RNN is able to express its uncertainty over the rate parameter of the demand Poisson process by computing a full posterior distribution. 
As qualitatively highlighted in Figure \ref{fig:prnn_vs_vprnn}, the VP-RNN predicts a full distribution over future demand rates possibly generating the observed data, whereas, by construction, the P-RNN only defines a point estimate for future rates.

\subsection{Prescriptive Results}
\label{subsec:inventory_management}
In this section, we focus on evaluating the inventory performance of different predictive models. Specifically, we count the number of shortages (of pickups and returns) during the next day based on actual demand data (i.e., the sequence consists of actual pickup and return events), assuming the initial bike inventory level prescribed by the UDF (see Section \ref{subsec:inventory_model}) obtained using the computed pickup and return forecasts. We set the unit penalty for each lost pickup and lost return to one ($l_p = l_r = 1$).

We compare our solutions to a model which receives perfect information about future demand patterns. 
Specifically, we use the UDF with perfect information about future rates $\hat{\mu}_t, \hat{\lambda}_t$.
Therefore, this approach serves as an \emph{oracle} that provides prescriptive performance in the limit of perfect forecasting accuracy for any algorithm within the same inventory decision model.

\begin{table}[t]
\caption{Test prescriptive performance for inventory management tasks averaged over all stations of interest.}
\small
\begin{tabular}{ccccccccc}
\hline
Aggregation & Measure & HA & MA & LR & P-RNN & VP-RNN & MOVP-RNN & Oracle \\
\hline
\multirow[t]{2}{*}{60 min} & Cost & 10.14 & 10.01 & 10.88 & 9.95 & \textbf{9.37} & 10.21 & 8.18 \\ 
& RPD & 24.1\% & 22.5\% & 33.1\% & 21.7\% & \textbf{14.6\%} & 24.9\% & - \\ 
\multirow[t]{2}{*}{30 min} & Cost & 10.13 & 10.01 & 10.88 & 9.56 & \textbf{8.91} & 10.25 & 8.18 \\ 
& RPD & 23.9\% & 22.4\% & 33.1\% & 16.8\% & \textbf{8.9\%} & 25.3\% & - \\
\multirow[t]{2}{*}{15 min} & Cost & 10.13 & 10.00 & 10.88 & 9.17 & \textbf{8.76} & 10.20 & 8.19 \\ 
& RPD & 23.7\% & 22.1\% & 32.8\% & 11.9\% & \textbf{7.0\%} & 24.5\% & - \\
\hline
\end{tabular}
\label{tab:decision_results}
\end{table}
In our experiments, we care about two main performance indicators: (i) \emph{Cost}: defined as the average number of unsatisfied customers per day due to the shortage of pickups and returns, (ii) \emph{Relative percentage difference} (RPD) from oracle performance, more formally:
$$
\mbox{RPD} = \frac{c_i - c_{oracle}}{c_{oracle}},
$$
where $c_{oracle}$, $c_i$ represent the costs obtained respectively by the oracle model and model $i \in $ \{HA, MA, LR, P-RNN, VP-RNN, MOVP-RNN\}.
Results in Table \ref{tab:decision_results} show that the forecasts generated by VP-RNN lead to the best decisions.
Specifically, inventory decisions based on VP-RNN predictions are able to decrease costs at least 40\% closer to oracle performance when compared to traditional HA, MA and LR performance and at least 10\% when compared to VP-RNN's deterministic counterpart P-RNN.
Results also highlight how predictive models which assume Poisson distributed demand (especially P-RNN and VP-RNN) benefit significantly from using smaller temporal discretizations. 
This is in line with the finding of \cite{Raviv&Kolka2013} that finer time discretizations yield a better fit of the non-homogeneous Poisson process to the pickup and return time series. Contrarily, classical models that do not use the Poisson property do not show a significant improvement when time intervals are chosen smaller.

Notably, our evaluation highlights a fundamental misalignment between prediction and decision performance.
Table \ref{tab:pred_results} shows that MOVP-RNN is able to obtain the best \emph{absolute} prediction performance across all stations and test days.
However, when evaluated in the context of decision performance, the predictions of MOVP-RNN lead to results close to the ones achieved by HA, MA and LR. 
\begin{figure}[t]
\centering
\includegraphics[width=\textwidth]{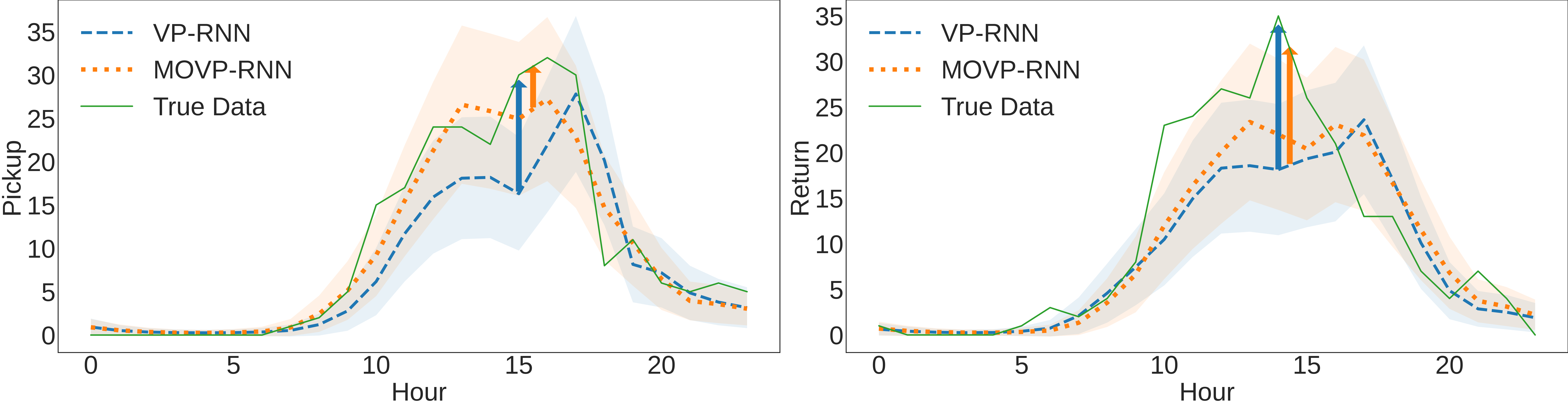}
\caption{A graphical representation of pickup (left) and return (right) predictions for one test day in station 3641 of the Citi Bike system. The plot shows how the two models reach different types of over/under-estimation patterns, where the arrows highlight the different biases during key moments of the day (e.g. the afternoon demand peak).}
\label{fig:so_vs_mo}
\end{figure}

\begin{table}[t]
\caption{Performance statistics for the predictions shown in Figure \ref{fig:so_vs_mo}. Results compare VP-RNN and MOVP-RNN on both prediction and decision performance.}
\small
\begin{tabular}{ccccccccc}
& \multicolumn{6}{c}{Prediction} & \multicolumn{2}{c}{Decision}\\
\hline 
& \multicolumn{3}{c}{Pickup} & \multicolumn{3}{c}{Return} &\\
& RMSE & MAE & $R^2$ & RMSE & MAE & $R^2$ & Inventory & Cost \\
\hline
VP-RNN & 5.35 & 3.53 & 0.76 & 6.02 & 3.92 & 0.69 & 11 & \textbf{13.72} \\
MOVP-RNN & \textbf{3.04} & \textbf{2.09} & \textbf{0.92} & \textbf{4.94} & \textbf{3.43} & \textbf{0.79} & 32 & 29.49 \\
\hline
\end{tabular}
\label{tab:so_vs_mo}
\end{table}
Figure \ref{fig:so_vs_mo} and Table \ref{tab:so_vs_mo} show a representative example where this misalignment is particularly evident.
Figure \ref{fig:so_vs_mo} compares pickup and return predictions belonging to VP-RNN and MOVP-RNN against real pickup and return observations on a single held-out test day, and Table \ref{tab:so_vs_mo} presents both prediction and decision performance for the same day.
In Figure \ref{fig:so_vs_mo}, we can observe how the two models have different error patterns, with the VP-RNN underestimating pickups and returns in a similar way, opposed to the MOVP-RNN which is approximately unbiased for the pickups but underestimates the return process.
From a strictly-predictive point of view, it is clear that MOVP-RNN represents a better model (Table \ref{tab:so_vs_mo}).
However, once the generated predictions are used in the UDF, the relative performance between the two models is reversed.

The reason behind this misalignment lies in the nature of the decision-making problem at hand. 
Specifically, when considering the task of selecting the best starting inventory, the optimal decision is fundamentally influenced by the cumulative difference between pickups and returns, rather than only their separate evolutions over the day. 
For example, in Figure \ref{fig:so_vs_mo}, by underestimating only the return rate, MOVP-RNN wrongly predicts a higher cumulative net demand (i.e., it predicts the correct number of pickups, but a lower number of returns), ultimately selecting higher starting inventories, and in practice leading to higher overall costs.
On the other hand, by having similar biases between pickup and return predictions, the VP-RNN will have a better estimate of the optimal starting inventory level, thus describing a situation where prediction errors in the same direction might (partially) cancel out when evaluated on decision performance.   

Current research fails to acknowledge these interactions between prediction and decision-making tasks, and rather focuses on either prediction or decision performance in isolation.
Motivated by this, we argue that having a deep understanding of how properties of the predictions can effect downstream decision-making processes is of fundamental importance. Therefore, we further explore this phenomenon in the remainder of this section.

\subsection{Qualitative analysis on the prediction-decision misalignment}
\label{subsec:qualitative_understanding}

To analyze the reason for the misalignment between predictive and prescriptive performance, we design a synthetic experiment. The goal of the experiment is to show how the bias of the prediction affects inventory decision making. We consider two types of bias: the \emph{same-side} bias, i.e. when pickups and returns are both either over-estimated or under-estimated, and the \emph{opposite-side} bias, i.e. when pickup and return estimates are biased in opposite directions. 
In what follows, we examine the impact of these different biasing patterns on the inventory decisions computed according to the model presented in Section \ref{subsec:inventory_model}. 

\begin{figure}[t]
\centering
\includegraphics[width=.6\textwidth]{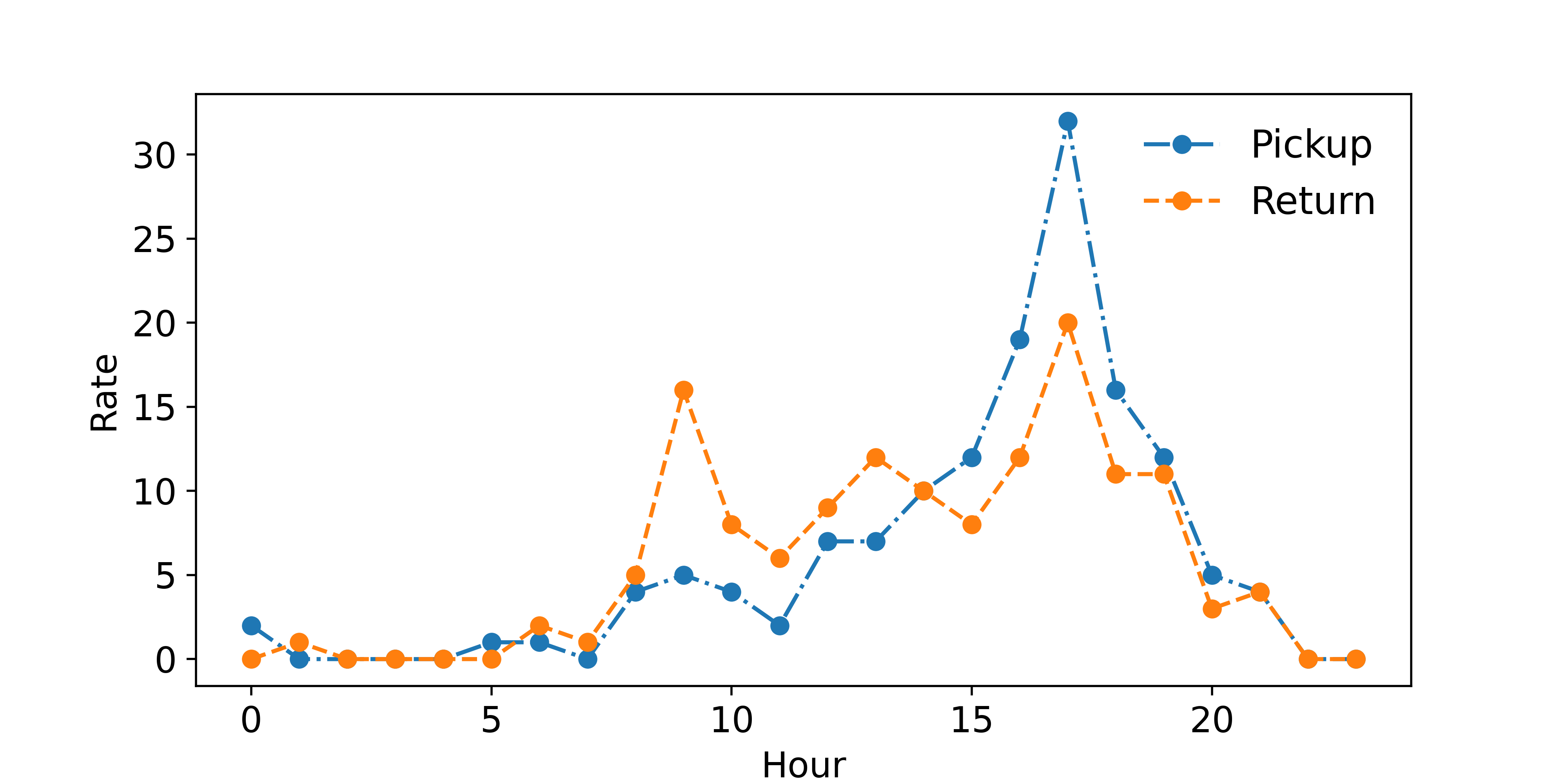}
\caption{Graphical representation of hourly pickup rates $\mu_t$ and return rates $\lambda_t$ for one test day of station 168 of the Citi Bike system, on which we base our synthetic experiment.}
\label{fig:demand}
\end{figure}

For exposition purposes, we select the demand pattern observed in station 168 on November 13, 2018. We choose a 60-minute aggregation time interval, and assume that we have perfect information about pickup and return patterns for that day. 
The observed counts of pickup and return arrivals within each hour $t$ are regarded as pickup rate $\mu_t$ and return rate $\lambda_t$, respectively. Figure \ref{fig:demand} shows the demand pattern of the selected instance, represented as pickup rate $\mu_t$ and return rate $\lambda_t, t \in[0,23]$. Given the information about the true count demand rates, the oracle starting inventory is then calculated according to \eqref{eq:trans_prob} - \eqref{eq:inventory_decision}.

\begin{figure}
\centering
\begin{subfigure}{.5\textwidth}
  \centering
  \includegraphics[width=\linewidth]{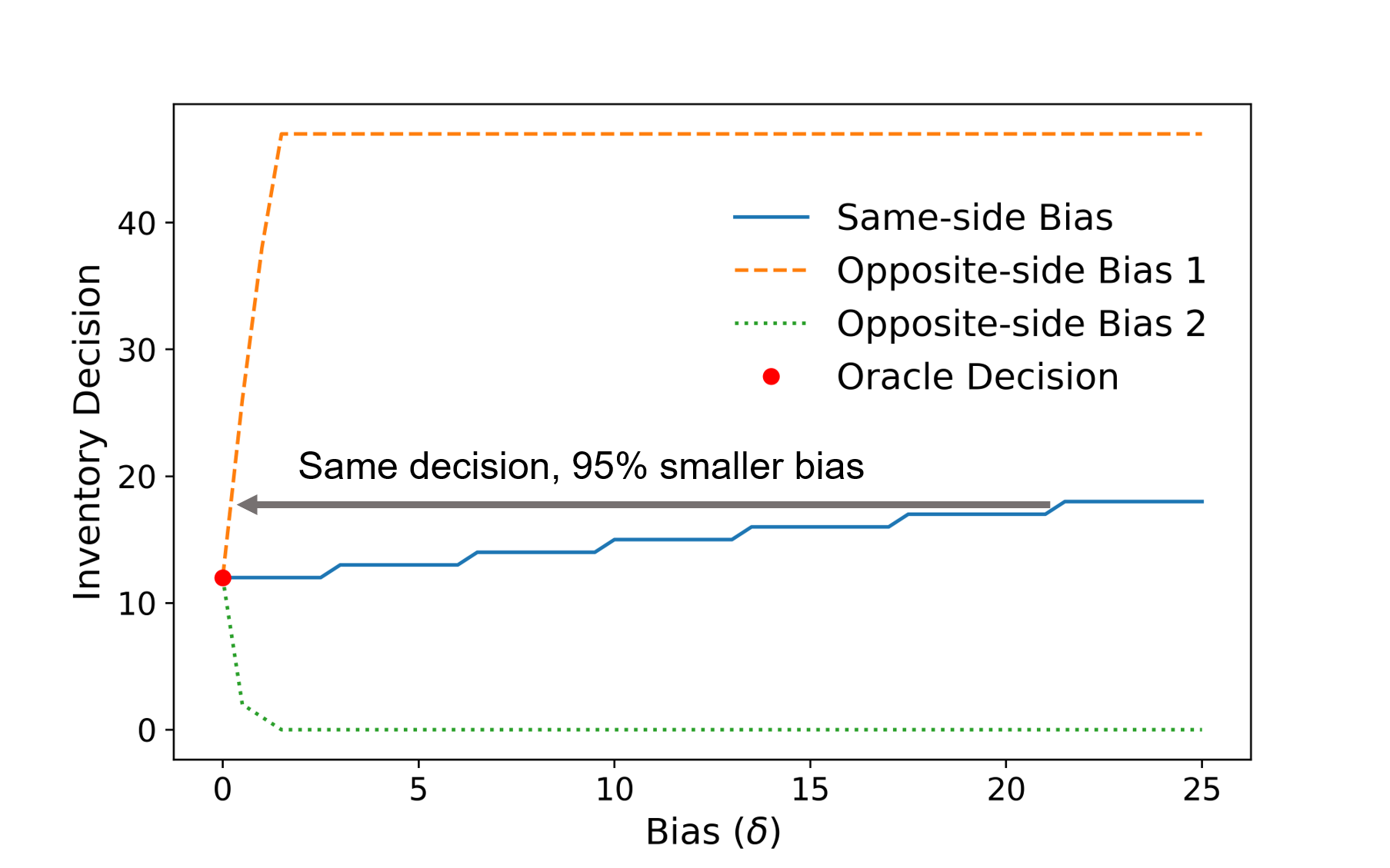}
  \caption{Inventory decision}
  \label{fig:sub1}
\end{subfigure}%
\begin{subfigure}{.5\textwidth}
  \centering
  \includegraphics[width=\linewidth]{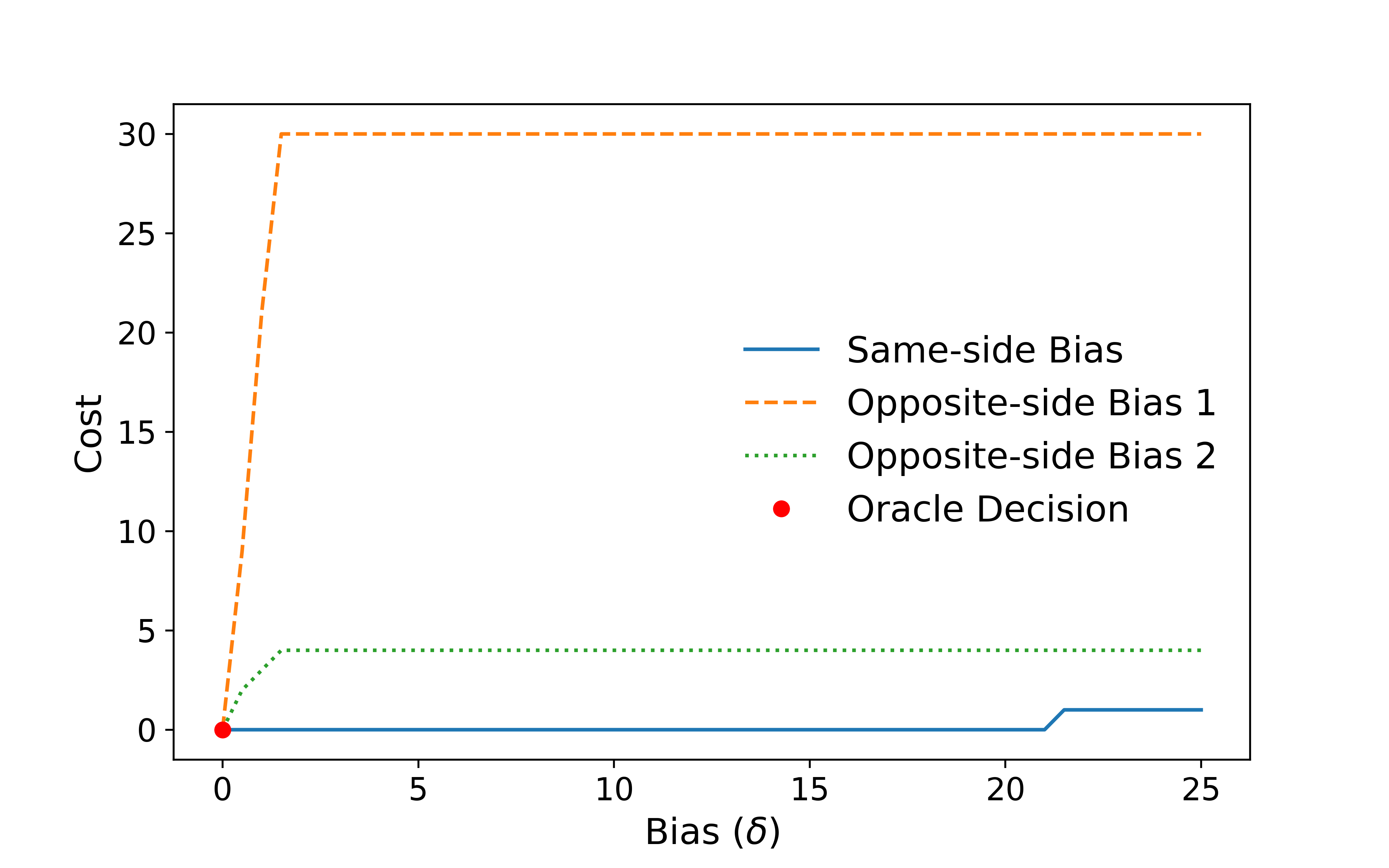}
  \caption{Decision performance}
  \label{fig:sub2}
\end{subfigure}
\caption{Impact of different bias patterns on inventory decisions as a function of bias level $\delta$. Plot (a) compares inventory decisions under bias predictions against the oracle inventory decision (red dot) under predictions with perfect information ($\delta=0$). To get the same deviation from the oracle inventory decision, the same-side bias can sustain a 95\% increase in forecasting error compared with the opposite-side bias. Plot (b) evaluates the performance of the decisions using the cost of lost sales.}
\label{fig:biased_demand}
\end{figure}

To assess the impact of different bias patterns, we select bias levels $\delta$ in the interval $[0,25]$ in increments of $0.5$. 
For each $\delta$ and each $t$, we generate the same-side biased pickup rate $\hat{\mu_t} = \mu_t + \delta$ and return rate $\hat{\lambda_t} = \lambda_t + \delta$. 
On the other hand, for the opposite-side bias, the pickup rate is over-estimated while the return rate is under-estimated by $\delta$, or vise versa. Specifically, the first opposite-side biased pickup rate and return rate are calculated as $\hat{\mu}_t = \mu_t + \delta$ and $\hat{\lambda}_t = \max(\lambda_t - \delta, 0)$. Note that in the case of downward biases, we truncate the resulting prediction at 0 to retain feasible estimates for the rates. Positive estimates are typically guaranteed by any prediction method. The second opposite-side biased pickup rate and return rate are calculated as $\hat{\mu}_t = \max(\mu_t - \delta, 0)$ and $\hat{\lambda}_t = \lambda_t + \delta$, respectively.
Finally, we calculate inventory decisions for all biased demand rates according to \eqref{eq:trans_prob} - \eqref{eq:inventory_decision}. 
Figure \ref{fig:biased_demand} illustrates the relationship between inventory decisions and bias level $\delta$, and the performance of decisions under bias level $\delta$. 

In Figure \ref{fig:biased_demand}, the inventory decision gradually deviates from the oracle decision as $\delta$ increases. Under the opposite-side bias 2, where pickup rates are under-estimated and return rates are over-estimated, the inventory decision quickly converges to the lower bound 0. Under the opposite-side bias 1, where pickups are over-estimated and returns are under-estimated, the inventory decision deviates rapidly towards the upper bound of station capacity. On the other hand, the inventory decision computed under the same-side bias is always very close to the oracle decision, even if the bias level $\delta$ is very large. Plot (b) shows that the resulting cost for any same-side bias smaller than 22 is 0, indicating that there are still no dissatisfied customers, despite the slightly different decision. Also for larger biases the cost only increases to 1. Contrarily, a small over-estimation of 2 units of pick-ups together with under-estimated returns already leads to 30 dissatisfied customers, a cost that is not caused by any same-side bias level in our experiment. Concluding, pickup and return rates with a very large bias on the same side still yield much better inventory decisions than those with a small bias on opposite sides. 

\subsection{Cumulative Error Metric}
\label{subsec:cumulative_error_metric}
The qualitative analysis highlights that improved forecasting accuracy on pickup and return rates individually, does not necessarily lead to better inventory decisions.
Specifically, the inventory decision model is based on a complex function of all pickup and return rates during the day and the difference between them, which creates a misalignment between prediction and decision objectives.
To deal with this misalignment, we introduce a new metric to measure forecast quality, the \emph{Cumulative Error} (CE) of net demand during a day.
We formally define CE as follows:

\begin{equation}\label{eq:sum_err}
    \text{CE} = \left\vert \sum_{t=0}^{T} [(\mu_t-\lambda_t) - (\hat{\mu_t}-\hat{\lambda_t})] \right\vert
\end{equation} 

In \eqref{eq:sum_err}, $\mu_t-\lambda_t$ represents the difference between the true pickup and return rates (i.e., actual net demand) within time interval $t$, while $\hat{\mu_t}-\hat{\lambda_t}$ represents the difference between the predicted pickup and return rates (i.e., predicted net demand) within time interval $t$. 

Table \ref{tab:sum_err} shows the average CE and the prescriptive decision performance over all test days and stations.
When evaluated on CE, VP-RNN clearly outperforms all other models, which is in line with its lowest corresponding costs. Even though its predictions of pickup and return rates are separately not the most accurate, VP-RNN still makes the best prescriptive decision.
This finding indicates that by measuring the cumulative error of net demand, CE yields a measure of prediction quality that is better aligned with the eventual decision performance.

\begin{table}[t]
\caption{Average CE and prescriptive performance over all test days and all 30 stations.}
\small
\begin{tabular}{cccccccc}
\hline
Aggregation & Measure & HA & MA & LR & P-RNN & VP-RNN & MOVP-RNN \\
\hline
\multirow[t]{3}{*}{60 min} & CE & 12.01 & 11.54 & 12.23 & 15.10 & \textbf{11.50} & 12.82\\
& Cost & 10.14 & 10.01 & 10.88 & 9.95 & \textbf{9.37} & 10.21\\ 
& RPD & 24.1\% & 22.5\% & 33.1\% & 21.7\% & \textbf{14.6\%} & 24.9\% \\  
\multirow[t]{3}{*}{30 min} & CE & 11.93 & 11.54 & 12.23 & 11.89 & \textbf{7.29} & 11.72\\
& Cost & 10.13 & 10.01 & 10.88 & 9.56 & \textbf{8.91} & 10.25\\ 
& RPD & 23.9\% & 22.4\% & 33.1\% & 16.8\% & \textbf{8.9\%} & 25.3\% \\
\multirow[t]{3}{*}{15 min} & CE & 11.93 & 11.54 & 12.23 & 11.16 & \textbf{6.19} & 10.95\\
& Cost & 10.13 & 10.00 & 10.88 & 9.17 & \textbf{8.76} & 10.20 \\ 
& RPD & 23.7\% & 22.1\% & 32.8\% & 11.9\% & \textbf{7.0\%} & 24.5\% \\
\hline
\end{tabular}
\label{tab:sum_err}
\end{table}

\section{Conclusion}
\label{sec:conclusion}
Taking effective operational-level decisions in real-world bike-sharing systems unavoidably entails a combination of demand forecasting with inventory decision making. 
Whereas both have separately received considerable attention, their interface has been left largely unaddressed. 
Bike pickups and returns jointly determine the inventory dynamics of a bike-sharing station, therefore, a classical forecast accuracy evaluation of both streams separately is not perfectly indicative for the quality of the resulting inventory decision. 
This paper illuminates this mismatch by considering a UDF to determine daily starting inventory levels in combination with various forecasting methods for the pickup and return rates. 
Among these are variations of a novel, deep generative model that represents pickup and return rates as latent variables, as well as several classical and learning-based benchmarks.

We show that the proposed method outperforms the benchmarks in terms of both forecast accuracy and the service quality of the resulting inventory decisions. 
Explicitly using a Poisson likelihood at prediction time and modelling the pickup and return rates as latent variables yields a better distributional fit and higher forecast accuracy on the 2018 Citi Bike data set, as evaluated using RMSE, MAE, and $R^2$. 
However, whereas these three measures agree on the model variant with the highest prediction accuracy (the MOVP-RNN), another variant of our approach (the VP-RNN) gives the best inventory decisions in terms of customer service quality. 
Classical measures of forecast accuracy, when used separately on pickup and return predictions, are not fully indicative of prescriptive (decision) accuracy. 
Our experiments show that conforming signs of the errors in pickups and returns during the day can to a large extent cancel out the consequences of their magnitudes.
Therefore, we propose a different accuracy measure, the Cumulative Error of net demand, and show that the ranking of forecasting methods based on this measure is in line with their ranking based on inventory performance.
We furthermore show -- in line with theoretical literature on non-homogeneous Poisson models for user dissatisfaction -- that our approach can additionally benefit from narrowing down the prediction interval to 30 or 15 minutes, whereas classical benchmarks cannot.

In the context of bike-sharing, this paper builds an intuition for what predictive properties are relevant for taking effective inventory decisions.
We highlight that 1) using simple averages, or even a linear regression, to estimate pickup and return rates leads to poor empirical decision performance, 2) learning-based approaches that exploit the Poisson likelihood in combination with a latent variable model for the rates lead to better decisions, but 3) it is crucial that predictive accuracy is measured in a way that aligns with the eventual inventory decision. 
That is, the joint effect of pickups and returns, their forecasts and biases should be taken into account.
This paper proposes an accuracy measure that accomplishes this in the often-encountered case of static rebalancing.

Future research should proceed on the interface of demand prediction and decision making for shared mobility. 
In the context of station-based bike-sharing, directly embedding decision-performance incentives within deep learning architectures would enable end-to-end learning of "decision-aware predictors". 
Integrating forecasting and decision-making further by means of imitation learning or data-driven optimization provides a different approach. 
Adapting the proposed accuracy measure to other objectives, such as dynamic rebalancing or joint, system-wide optimization will lead to new insights on forecasting methods' performances in these scenarios.
Efficiently exploiting decision performance in system-level forecasting of demand patterns, auto- and cross-sectional correlations, remains an open challenge.
Possible extensions to free-floating or hybrid systems would open a plethora of new application areas, such as cars, scooters, and urban air mobility.

\section*{Acknowledgements}
This research was supported by TUM International Graduate School of Science and Engineering (IGSSE) through the project ILOMYTS. Dennis Prak was supported by the Dutch Research Council (NWO) [grant nr. 019.191SG.003].

\bibliographystyle{apalike}
\bibliography{ref}

\appendix
\section{ELBO derivation of \eqref{eq:elbo}}
\label{appendix:elbo}
We hereby report the derivation of the evidence lower bound used for inference in the VP-RNN: 
\begin{align*}
    \log p_{\theta}(\mathbf{x}_{1:T}) & = \log \int p_{\theta}(\bx_{1:T}, \boldsymbol{\lambda}_{1:T}, \bh_{1:T}) \mathrm{d}\lambda \, \mathrm{d}\bh \\
    & = \log \int \frac{q_{\phi}(\boldsymbol{\lambda}_{1:T} \, | \, \bx_{1:T})}{q_{\phi}(\boldsymbol{\lambda}_{1:T} \, | \, \bx_{1:T})} p_{\theta}(\bx_{1:T}, \boldsymbol{\lambda}_{1:T}, \bh^p_{1:T}) \mathrm{d}\lambda \, \mathrm{d}\bh\\
    & = \log \mathbb{E}_{q_{\phi}(\boldsymbol{\lambda}_{1:T} \, | \, \bx_{1:T})} \left[\prod_{t=1}^{T} \frac{p_{\theta}(\bx_{t} \, | \, \boldsymbol{\lambda}_{t}) p_{\theta}(\boldsymbol{\lambda}_{t} \, | \, \bh^p_{t}) p_{\theta}(\bh^p_{t} \, | \, \bh^p_{t-1}, \mathbf{u}_{t})}{q_{\phi}(\boldsymbol{\lambda}_{t} \, | \, \bh^q_{t}, \bx_{t}, \bu_t)} \right] \\
    & \geq \mathbb{E}_{q_{\phi}(\boldsymbol{\lambda}_{1:T} \, | \, \bx_{1:T})} \left[\sum_{t=1}^{T} \log p_{\theta}(\bx_{t} \, | \, \boldsymbol{\lambda}_{t}) + \log p_{\theta}(\bh^p_{t} \, | \, \bh^p_{t-1}, \mathbf{u}_{t}) + \log \left( \frac{p_{\theta}(\boldsymbol{\lambda}_{t} \, | \, \bh^p_{t})}{q_{\phi}(\boldsymbol{\lambda}_{t} \, | \, \bh^q_{t}, \bx_{t}, \bu_t)} \right) \right] \\
    & = \mathbb{E}_{q_{\phi}(\boldsymbol{\lambda}_{1:T} \, | \, \bx_{1:T})} \left[\sum_{t=1}^{T} \log p_{\theta}(\bx_{t} \, | \, \boldsymbol{\lambda}_{t}) + \log p_{\theta}(\bh^p_{t} \, | \, \bh^p_{t-1}, \mathbf{u}_{t}) \right] \\
    & - \sum_{t=1}^{T} \mathbb{KL} \left(q_{\phi}(\boldsymbol{\lambda}_{t} \, | \, \bh^q_{t}, \bx_{t}, \bu_t) \,\mid\mid\, p_{\theta}(\boldsymbol{\lambda}_{t} \, | \, \bh^p_{t}) \right) \\ 
\end{align*}

\end{document}